\documentclass[12pt,leqno]{amsart}
\usepackage{amssymb}
\usepackage{amsmath,amssymb,color}
\numberwithin{equation}{section}

\newcommand{\OOO}{\Omega}

\newcommand{\la}{\lambda}
\newcommand{\va}{\varphi}
\newcommand{\ppp}{\partial}

\newcommand{\sumk}{\sum_{k=1}^{\infty}}

\newcommand{\EEE}{E_{\alpha,1}}

\newcommand{\R}{\mathbb{R}}
\newcommand{\C}{\mathbb{C}} 
\newcommand{\N}{\mathbb{N}}

\newcommand{\ooo}{\overline}

\newcommand{\pppa}{\ppp_t^{\alpha}}

\allowdisplaybreaks

\title
[]
{
Uniqueness in determining the orders of time and spatial
fractional derivatives
}

\author{by Masahiro Yamamoto}

\thanks{
Graduate School of Mathematical Sciences, The University
of Tokyo, Komaba, Meguro, Tokyo 153-8914, Japan \\
Honorary Member of Academy of Romanian Scientists, 
Splaiul Independentei Street, no 54,
050094 Bucharest Romania \\
Peoples' Friendship University of Russia 
(RUDN University) 6 Miklukho-Maklaya St, Moscow, 117198, Russian Federation
e-mail: {\tt myama@ms.u-tokyo.ac.jp}\\
}

\date{}
\begin{document}
\maketitle
\baselineskip 18pt
\begin{abstract}
We prove the uniqueness in determining both orders of 
fractional time derivatives and spatial derivatives in diffusion equations 
by pointwise data.  The proof relies on the eigenfunction expansion and 
the asymptotics of the Mittag-Leffler function. 
\end{abstract}

\baselineskip 18pt

\noindent

\section{Introduction and main result}
 
Let $\OOO \subset \R^d$ be a bounded domain with smooth boundary, and 
let $\alpha \in (0,2) \setminus \{ 1\}$ and $0<\beta < 1$.
For suitably given $a \in L^2(\OOO)$, let $u_{\alpha,\beta}
= u_{\alpha,\beta}(x,t)$ satisfy 
$$\left\{ \begin{array}{rl}
& \pppa u = -A^{\beta} u, \quad x \in \OOO, \, t>0,\\
& u\vert_{\ppp\OOO} = 0, \qquad t>0, \\
& \left\{\begin{array}{rl}
  &u(x,0) = a(x), \quad x \in \OOO \quad \mbox{if $0<\alpha<1$},\\
  & u(x,0) = a(x), \quad \ppp_tu(x,0) = 0 \quad \mbox{if $1<\alpha<2$}.
  \end{array}\right.                             
\end{array}\right.                     \eqno{(1.1)}
$$
Here $\ppp_t^{\alpha}$ denotes the Caputo derivative:
$$
\ppp_t^{\alpha}v(t) = \frac{1}{\Gamma(n-\alpha)}\int^t_0 (t-s)^{n-\alpha-1}
\frac{d^nv}{ds^n}(s) ds \quad \mbox{for $n-1<\alpha<n$ with $n\in \N$}.
$$
Moreover we set 
$$
Au(x) = -\sum_{i,j=1}^d \ppp_i(a_{ij}(x)\ppp_ju(x)) - c(x)u(x), \quad 
x\in \OOO
$$
where $a_{ij} = a_{ji} \in C^1(\ooo{\OOO})$ for $1\le i,j\le n$ and
$c \in C(\ooo{\OOO})$, $\le 0$ in $\OOO$, and there exists a constant 
$\mu>0$ such that 
$$
\sum_{i,j=1}^d a_{ij}(x)\xi_i\xi_j \ge \mu\sum_{i=1}^d \xi_i^2,
\quad x\in \ooo{\OOO}, \, \xi_1, ...., \xi_d \in \R.
$$
To $A$, we attach the domain of the operator:
$$
D(A) = H^2(\OOO) \cap H^1_0(\OOO),
$$
where $H^2(\OOO)$ and $H^1_0(\OOO)$ are usual Sobolev spaces.

By the symmetry and positivity of $A$, we can 
define the fractional power $A^{\beta}$ of $A$ (e.g., Pazy \cite{Pa},
Tanabe \cite{Ta}).  We can interpret $2\beta$ as the order of
the spatial fractional derivative because the differential oprerator $A$ 
is of second order.
A special case $\beta=1$ in (1.1) describes a time-fractional diffusion 
equation, and is a macroscopic model for the continuous-time random walk, 
which is a meaningful model e.g., for anomalous diffusion in heterogeneous 
media.  As for phyiscal backgrounds, see Metzler and Klafter \cite{MK} for 
example.

There has been the rapidly increasing literature on mathematical 
works for initial boundary value problem (1.1) for time fractional diffusion 
equations and we refer only to a very limited number of articles:
Gorenflo, Luchko and Yamamoto \cite{GLY}, Kubica, Ryszewska and Yamamoto
\cite{KRY}, Kubica and Yamamoto \cite{KY}, Sakamoto and Yamamoto
\cite{SY}, Zacher \cite{Za}.

For analyzing such an initial boundary value problem (1.1), we need to 
discuss inverse problems which are concerned with a quantitative 
procedure for the determination of parameters in (1.1) such as 
$\alpha$, $\beta$, $a_{ij}(x)$, etc.

In particular, the orders $\alpha$ and $\beta$ are essential for modelling 
the phenomena under consideration.
The main purpose of this article is to establish the uniqueness for 
\\
{\bf Inverse problem.}\\
{\it Let $x_0 \in \OOO$ be given.  Determine $\alpha \in (0,2) \setminus 
\{1\}$ and $\beta \in (0,1)$ by $u_{\alpha,\beta}(x_0,t)$, $0<t<T$.
}
\\

For the statement of our main result, we need to introduce the eigensystem of 
the elliptic operator $A$ of the second order.
Let 
$$
0<\la_1 < \la_2 < .... \longrightarrow \infty
$$
be the set of all the eigenvalues of $A$, and 
let $\va_{kj}$, $1\le j\le m_k$ be an orthonormal system of eigenfunctions of 
$A$ for $\la_k$, $k\in \N$: $A\va_{kj} = \la_k\va_{kj}$, $1\le j \le m_k$.

We further define the Mittag-Leffler function by 
$$
E_{\alpha,1}(z) = \sum_{k=0}^{\infty} \frac{z^k}{\Gamma(\alpha k + 1)},
\quad z \in \C.
$$
It is known that $E_{\alpha,1}(z)$ is an entire function in $z \in \C$
(e.g., Podlubny \cite{Po}).

Now we are ready to state our main result.
\\
Thus we have established 
\\
{\bf Theorem (uniqueness)}.\\
{\it Let $\gamma > \frac{d}{4}$.  We assume that 
$$
a\ge 0, \not\equiv 0 \quad \mbox{or} \quad \le 0, \not\equiv 0 \quad 
\mbox{in $\OOO$},
\quad a \in H^{2\gamma}_0(\OOO)            \eqno{(1.2)}
$$
and there exists $k_0\in \N$ such that 
$$
\sum_{j=1}^{m_{k_0}} (a, \va_{k_0j})\va_{k_0j}(x_0) \ne 0, \quad \la_{k_0} \ne 1.                                                            \eqno{(1.3)}
$$
Then $u_{\alpha,\beta}(x_0,t)$ for $0<t<T$ uniquely determines
$\alpha \in (0,2)\setminus \{ 1\}$ and $\beta \in (0,1)$.
}
\\

Condition (1.3) is essential for the uniqueness of $\beta$.
Indeed, for simplicity, let each $\la_k$ be a simple eigenvalue:
$m_k=1$ for $k \in \N$.  We write $\va_k := \va_{km_k} = \va_{k1}$
for $k \in \N$.  Let $\la_1 = 1$ and $a=\va_1$.  Then $A^{\beta}a
= a$ for any $\beta>0$.
Then we can readily verify that $v(x,t) = E_{\alpha,1}(-t^{\alpha})
\va_1(x)$ satisfies (1.1) with $a = \va_1$: $u_{\alpha,\beta}(x,t) 
= E_{\alpha,1}(-t^{\alpha})\va_1(x)$.
Therefore we cannot determine $\beta$ by $u_{\alpha,\beta}(x_0,t) 
= E_{\alpha,1}(-t^{\alpha})\va_1(x_0)$, because the solution is 
independent of $\beta$.

The condition (1.3) is satisfied if $a(x_0) \ne 0$ and there exists an 
eigenvalue not $1$.

Indeed we have  
$$
a(x_0) = \sum_{k=1}^{\infty} \sum_{j=1}^{m_k} (a,\va_{kj})\va_{kj}(x_0).
$$
If $\sum_{j=1}^{m_k} (a,\va_{kj})\va_{kj}(x_0) = 0$ for all $k\in \N$, then 
$a(x_0) = 0$.  

For the uniqueness for $\alpha$, as is seen in the proof, we need not
(1.2) itself and the uniqueness follows from $(A^{-\beta}a)(x_0) \ne 0$.
\\

As for inverse problems of determining order $\alpha$, we 
explain existing related works.  
Hatano, Nakagawa, Wang and Yamamoto \cite{HNWY} is an early theoretical 
work and we further refer to Ashurov and Umarov \cite{AU}, Janno \cite{J},
Janno and Kinash \cite{JK}, Krasnoschok, Pereverzyev, Siryk and Vasylyeva
\cite{KPSV}, Li and Yamamoto \cite{LY}, 
Yu, Jing and Qi \cite{YJQ}.  Moreover see a survey chapter Li, Liu and 
Yamamoto \cite{LiLiuY}.  
In the case of $\OOO = (0,1)$ and $A = -\Delta$ in (1.1), 
the article Tatar and Ulusoy \cite{TU}
proves the uniqueness in determining $\alpha$ and $\beta$
simultaneously by $u(x_0,t)$, $0<t<T$ with fixed $x_0 \in (0,1)$
under the assumption 
$$
(a,\va_k)_{L^2(0,1)} > 0 \quad \mbox{for all $k\in \N$}.   
                                    \eqno{(1.4)}
$$
Here we note that in the one-dimensional $\OOO$, all the eigenvalues 
are simple: $m_k=1$ and we write $\va_{km_k} = \va_k$.
The condition (1.4) requires the non-vansihing of all the Fourier coefficients
and is essentially different from our conditions (1.2) - (1.3), 
\\

In addition to the inverse problems of determining orders, there are 
very many works on other types of inverse problems for fractional equations, 
and we here refer to very limited related articles:
Cheng, Nakagawa, Yamamoto and Yamazaki \cite{CNYY},
Li, Zhang, Jia and Yamamoto \cite{LZJY}, Li, Luchko and Yamamoto \cite{LLY}, 
Ruan, Zhang and Wang \cite{RZW},
Rundell and Zhang \cite{RZ}, 
Sun, Li and Jia \cite{SLJ}, Zhang, Cheng and Wang \cite{ZCW},
two survey papers Li and Yamamoto \cite{LiY}, Liu, Li and Yamamoto
\cite{LiuLiY}.

This article is composed of three sections.  In Section 2, we will prove the 
theorem.
Section 3 is devoted to concluding remarks.

\section{Proof of Theorem}

{\bf First Step.}\\
In terms of the eigenvalues and the eigenfunctions, we can represent
$$
u_{\alpha,\beta}(x,t) = \sumk E_{\alpha,1}(-\la_k^{\beta}t^{\alpha})
\sum_{j=1}^{m_k} (a,\va_{kj})\va_{kj}(x), \quad x\in \OOO, \, 
t>0
$$
pointwise under the regularity condition of $a$ in the theorem.
More precisely, we have
$$
\Vert A^{\gamma}u_{\alpha,\beta}(\cdot,t)\Vert_{L^2(\OOO)}
\le C\Vert A^{\gamma}a\Vert_{L^2(\OOO)}
$$ 
and so we see that 
$u_{\alpha,\beta}(\cdot,t) \in C(\ooo{\OOO})$ by the 
Sobolev embedding and $4\gamma > d$.

We assume that
$$
u_{\alpha,\beta}(x_0,t) = u_{\alpha_1,\beta_1}(x_0,t), \quad 
0<t<T,
$$
for $\alpha, \alpha_1 \in (0,2) \setminus \{1\}$ and 
$\beta, \beta_1 \in (0,1)$.
Substituting $x=x_0$, we assume 
$$
\sumk E_{\alpha,1}(-\la_k^{\beta}t^{\alpha})
\sum_{j=1}^{m_k} (a,\va_{kj})\va_{kj}(x_0)
= \sumk E_{\alpha_1,1}(-\la_k^{\beta_1}t^{\alpha_1})
\sum_{j=1}^{m_k} (a,\va_{kj})\va_{kj}(x_0), \quad t>0.
$$
By the asymptotics of $\EEE(-z)$ with large $z>0$
(Theorem 1.4: pp.33-34 in Podlubny \cite{Po}), we obtain
$$
u_{\alpha,\beta}(x,t) 
= \sumk \sum_{j=1}^{m_k} \frac{(a,\va_{kj})\va_{kj}(x)}
{\Gamma(1-\alpha)\la_k^{\beta}t^{\alpha}}
+ O\left( \frac{1}{t^{2\alpha}}\right)
$$
$$
= \frac{1}{\Gamma(1-\alpha)t^{\alpha}}(A^{-\beta}a)(x) 
+ O\left( \frac{1}{t^{2\alpha}}\right) \quad \mbox{as $t \to \infty$.}
                                                                \eqno{(2.1)}
$$
We note that the fractional power $A^{-\beta}$ is defined by 
$$
A^{-\beta}a = \frac{\sin \pi\beta}{\pi}
\int^{\infty}_0 \eta^{-\beta}(A+\eta)^{-1}a d\eta
\quad \mbox{in $L^2(\OOO)$}, \quad 0<\beta<1                  \eqno{(2.2)}
$$
(e.g., Pazy \cite{Pa}, Tanabe \cite{Ta}).

{\bf Second Step.}\\
We here prove
\\
{\bf Lemma.}\\
{\it Let $a\ge 0$, $\not\equiv 0$ in $\OOO$ and $a \in H^{2\gamma}_0(\OOO)$
with $\gamma > \frac{d}{4}$.
Then 
$$
(A^{-\beta}a)(x) > 0, \quad x\in \OOO.
$$
}
\\
{\bf Proof.}\\
By the strong maximum principle for $\Delta - \eta$ with $\eta\ge 0$, we see
$$
(A+\eta)^{-1}a(x) >  0, \quad x \in \OOO, \, \eta \ge 0.  \eqno{(2.3)}
$$
Indeed we set $w = (A+\eta)^{-1}a$, that is,
$$
\left\{ \begin{array}{rl}
& \sum_{i,j=1}^d \ppp_i(a_{ij}(x)\ppp_jw) - \eta w = -a \le 0 
\quad \mbox{in $\OOO$}, \\
& w\vert_{\ppp\OOO} = 0.
\end{array}\right.
$$
The maximum principle (Corollary 4.3 in Renardy and Rogers \cite{RR},
p.106) yields 
$$
\min_{\ooo{\OOO}} w \ge \min_{\ppp\OOO} \min\{ w, 0\} = 0
$$
and so
$$
w \ge 0 \quad \mbox{on $\ooo{\OOO}$}.             \eqno{(2.4)}
$$
The strong maximum principle (\cite{RR}, p.109) implies that 
$$
\mbox{$w$ cannot achieve a non-positive minimum at any point in $\OOO$}, 
                                               \eqno{(2.5)}
$$
if $w$ is not a constant function.

Since $w\vert_{\ppp\OOO} = 0$, if $w$ is constant, then $w\equiv 0$, which is
impossible by $a \not\equiv 0$.  Therefore $w$ is not constant.
Therefore (2.5) holds.

By (2.4) and $w\vert_{\ppp\OOO} = 0$, we see that $\min_{\ooo{\OOO}} w = 0$.

Assume that there exists $x_0 \in \OOO$ such that $w(x_0) = 0$, that is,
$w$ attains a non-positive minimum (i.e., $0$) at $x_0\in \OOO$. 
This is impossible by (2.5).  Thus (2.3) is verified.

By $a \in D(A^{\gamma}) \subset H^{2\gamma}(\OOO) \subset C(\ooo{\OOO})$ 
with $\gamma > \frac{d}{4}$,  by (2.2) we have 
\begin{align*}
& A^{\gamma}\left(A^{-\beta}a 
- \frac{\sin \pi\beta}{\pi}
\int^{\infty}_0 \eta^{-\beta}(A+\eta)^{-1}a d\eta\right) \\
=& A^{-\beta}(A^{\gamma}a) 
- \frac{\sin \pi\beta}{\pi}
\int^{\infty}_0 \eta^{-\beta}(A+\eta)^{-1}(A^{\gamma}a) d\eta,
\end{align*}
and so (2.2) holds in $C(\ooo{\OOO})$.  We can substitute 
$x=x_0$ in (2.2) to conclude that $(A^{-\beta}a)(x)>0$ for 
all $x \in \OOO$.  Thus the proof of Lemma is complete.

{\bf Third Step.}\\
Since $u_{\alpha,\beta}(x_0,t) = u_{\alpha_1,\beta_1}(x_0,t)$ for 
$0 < t < T$, the $t$-analyticity means 
$u_{\alpha,\beta}(x_0,t) = u_{\alpha_1,\beta_1}(x_0,t)$ for $t>0$.
It follows from (2.1) that 
$$
\frac{1}{\Gamma(1-\alpha)t^{\alpha}}(A^{-\beta}a)(x_0) 
+ O\left( \frac{1}{t^{2\alpha}} \right)
= \frac{1}{\Gamma(1-\alpha_1)t^{\alpha_1}}(A^{-\beta_1}a)(x_0) 
+ O\left( \frac{1}{t^{2\alpha_1}}\right)
$$
as $t \to \infty$.  By Lemma, we see that $p := (A^{-\beta}a)(x_0) \ne 0$ and
$p_1 := (A^{-\beta_1}a)(x_0) \ne 0$.
Let $\alpha < \alpha_1$.  Then
$$
\frac{p}{\Gamma(1-\alpha)} + O\left( \frac{1}{t^{\alpha}}\right)
= \frac{p_1}{\Gamma(1-\alpha_1)t^{\alpha_1-\alpha}} 
+ O\left( \frac{1}{t^{2\alpha_1-\alpha}}\right)
$$
as $t \to \infty$.  Letting $t \to \infty$, by $\alpha_1>\alpha$, we 
reach $\frac{p}{\Gamma(1-\alpha)} = 0$, which is a contradiction
against $p\ne 0$.
Therefore $\alpha_1 \le \alpha$.  Similarly we can prove $\alpha_1
\ge \alpha$.  Thus we can conclude $\alpha_1 = \alpha$.

{\bf Fourth Step.}\\
Finally we will prove $\beta = \beta_1$.  
Since 
$$
u_{\alpha,\beta}(x_0,t) = \sum_{k=1}^{\infty}\sum_{j=1}^{m_k}
(a,\va_{kj})\EEE(-\la_k^{\beta}t^{\alpha})\va_{kj}(x_0), \quad 
0<t<T,
$$
is analytic in $t$ and convergent pointwise, again by the asymptotics
(e.g., Theorem 1.4 (pp.34-35) in \cite{Po}), for any $N\in \N$, we have
\begin{align*}
& \sum_{\ell=1}^N \frac{(-1)^{\ell+1}}{\Gamma(1-\alpha\ell)}
\frac{1}{t^{\alpha\ell}} \sum_{k=1}^{\infty} \frac{1}{\la_k^{\beta\ell}} 
\sum_{j=1}^{m_k} (a,\va_{jk})\va_{jk}(x_0)
+ O\left( \frac{1}{t^{\alpha(N+1)}}\right) \\
= & \sum_{\ell=1}^N \frac{(-1)^{\ell+1}}{\Gamma(1-\alpha\ell)}
\frac{1}{t^{\alpha\ell}} \sum_{k=1}^{\infty} \frac{1}{\la_k^{\beta_1\ell}} 
\sum_{j=1}^{m_k} (a,\va_{jk})\va_{jk}(x_0)
+ O\left( \frac{1}{t^{\alpha(N+1)}}\right)
\end{align*}
as $t \to \infty$.
For short descriptions, we set $a_k=
\sum_{j=1}^{m_k} (a,\va_{jk})\va_{jk}(x_0)$, $k\in \N$.
Then
\begin{align*}
& \sum_{\ell=1}^N \frac{(-1)^{\ell+1}}{\Gamma(1-\alpha\ell)}
\frac{1}{t^{\alpha\ell}} \sum_{k=1}^{\infty}
\frac{a_k}{\la_k^{\beta\ell}} + O\left( \frac{1}{t^{\alpha(N+1)}}\right)\\
=& \sum_{\ell=1}^N \frac{(-1)^{\ell+1}}{\Gamma(1-\alpha\ell)}
\frac{1}{t^{\alpha\ell}} \sum_{k=1}^{\infty}
\frac{a_k}{\la_k^{\beta_1\ell}} + O\left( \frac{1}{t^{\alpha(N+1)}}\right)
\end{align*}
as $t \to \infty$ for all $N\in \N$.
Since $\frac{1}{\Gamma(1-\alpha\ell)} > 0$ if 
$1-\alpha\ell \not\in \{ 0,-1,-2,-3, ....\}$, we obtain
$$
\sum_{k=1}^{\infty} \frac{a_k}{\la_k^{\beta\ell}}
= \sum_{k=1}^{\infty} \frac{a_k}{\la_k^{\beta_1\ell}}, \quad 
\ell \in \N \setminus \left\{ \frac{m}{\alpha} \right\}
_{m\in \N}.                       \eqno{(2.6)}
$$

Let $k_0$ be the minimum natural number satisfying (1.3):
$a_{k_0} \ne 0$.
Then 
$$
\sum_{k=k_0}^{\infty} \frac{a_k}{\la_k^{\beta\ell}}
= \sum_{k=k_0}^{\infty} \frac{a_k}{\la_k^{\beta_1\ell}}, \quad 
\ell \in \N \setminus \left\{ \frac{m}{\alpha} \right\}
_{m\in \N}. 
$$
\\
{\bf Case 1: $\la_{k_0} > 1$.}\\
We assume that $\beta_1 > \beta$.
Multiplying (2.6) with $\la_{k_0}^{\ell\beta}$, we have
$$
a_{k_0}( 1 - (\la_{k_0}^{\beta-\beta_1})^{\ell})
+ \sum_{k=k_0+1}^{\infty} a_k \left( \left( \frac{\la_{k_0}^{\beta}}
{\la_k^{\beta}} \right)^{\ell}
- \left( \frac{\la_{k_0}^{\beta}}{\la_k^{\beta_1}} \right)^{\ell}\right)
= 0                                    \eqno{(2.7)}
$$
for $\ell \in \N \setminus \left\{ \frac{m}{\alpha} \right\}
_{m\in \N}$.   
By $\beta - \beta_1 < 0$, $\la_{k_0} > 1$ and 
$\left\vert \frac{\la_{k_0}}{\la_k}\right\vert < 1$ for $k \ge k_0+1$, 
we have
$\vert \la_{k_0}^{\beta-\beta_1} \vert < 1$ for $k\in \N$ and 
$\left\vert \left( \frac{\la_{k_0}}{\la_k} \right)^{\beta} \right\vert$,
$\left\vert \frac{\la_{k_0}^{\beta}}{\la_k^{\beta_1}} \right\vert < 1$ for
$k\ge k_0+1$. Then we choose a sequence $\ell_n 
\in \N \setminus \left\{ \frac{m}{\alpha} \right\}_{m\in \N}$ 
such that $\lim_{n\to\infty} \ell_n = \infty$.

Indeed we assume that such a sequence $\ell_n$ does not exist.
Then there exists $N_0 \in \N$ such that 
$$
\left\{ \frac{m}{\alpha} \right\}_{m\in \N} \supset \{N_0, N_0+1, N_0+2,
... \}.
$$  
Therefore we can choose $m, m' \in \N$ such that $\frac{m}{\alpha}
= N_0$ and $\frac{m'}{\alpha} = N_0+1$.  Hence the subtraction 
yields $\frac{m'-m}{\alpha}= 1$, that is, $\alpha=m'-m$.
This is a contradiction because $0<\alpha<2$ and $\alpha \ne 1$.

Then, applying the 
Lebesgue convergence theorem for the series 
to the second term on the left-hand side
of (2.7), we reach $a_{k_0} = 0$.  This is impossible and so 
$\beta_1 \le \beta$.  Similarly we can prove $\beta_1 \ge \beta$.
Therefore we proved $\beta_1 = \beta$ in the case of $\la_{k_0} > 1$.
\\
{\bf Case 2: $\la_{k_0} < 1$.}

We assume that $\beta_1 > \beta$.
we multiplying (2.6) with $\la_{k_0}^{\ell\beta_1}$.
Then 
$$
a_{k_0}\left( \left( \frac{\la_{k_0}^{\beta_1}}{\la_{k_0}^{\beta}}
\right)^{\ell}
- 1 \right) 
+ \sum_{k=k_0+1}^{\infty} a_k \left( \left( \frac{\la_{k_0}^{\beta_1}}
{\la_k^{\beta}} \right)^{\ell}
- \left( \frac{\la_{k_0}^{\beta_1}}{\la_k^{\beta_1}} \right)^{\ell}\right)
= 0                            
$$
for $\ell \in \N \setminus \left\{ \frac{m}{\alpha} \right\}
_{m\in \N}$.   
Since $\beta_1 - \beta > 0$, $\la_{k_0} < 1$ and $\la_k > \la_{k_0}$ for
$k \ge k_0+1$, we have
$$
\left\vert \frac{\la_{k_0}^{\beta_1}}{\la_{k_0}^{\beta}}\right\vert < 1,
\quad \left\vert \frac{\la_{k_0}^{\beta_1}}{\la_k^{\beta}} \right\vert,\,
\left\vert \frac{\la_{k_0}^{\beta_1}}{\la_k^{\beta_1}} \right\vert < 1
$$
for $k \ge k_0+1$.  Letting $\ell \to \infty$ within 
$\N \setminus \left\{ \frac{m}{\alpha} \right\}_{m\in \N}$, 
similarly to Case 1, we see that $-a_{k_0} = 0$, 
which is impossible.  Therefore $\beta_1 \le \beta$.
Similarly we can prove $\beta_1 \ge \beta$.
Therefore we proved $\beta_1 = \beta$ in both Cases 1 and 2.
Thus the proof of Theorem is complete.

\section{Concluding remarks.}
\begin{itemize}
\item
We can discuss the inverse problem by data
$$
\int_{\OOO} u_{\alpha,\beta}(x,t) \rho(x) dx, \quad 0<t<T
$$
with suitable weight $\rho(x)$
or 
$$
\frac{\ppp u_{\alpha,\beta}}{\ppp \nu_A}(x_0,t), \quad 0<t<T,
$$
where $x_0 \in \ppp\OOO$ is suitably given and 
$\frac{\ppp v}{\ppp \nu_A}$ denotes the conormal derivative:
$\frac{\ppp v}{\ppp \nu_A}(x) = \sum_{i,j=1}^d a_{ij}(x)\nu_j(x)$
with the outward unit normal vector $(\nu_1(x), ..., \nu_d(x))$ to 
$\ppp\OOO$.
\item
The uniqueness is the primary theoretical issue for the inverse problem.
The stability is also important but it is open so far.
\end{itemize}

\section*{Acknowledgments}
The author thanks Professor Takeshi Saito (The University of Tokyo)
for the improvement of the argument in the proof of Theorem. 
The author was supported by Grant-in-Aid for Scientific Research (S)
15H05740 of Japan Society for the Promotion of Science and
by The National Natural Science Foundation of China
(no. 11771270, 91730303).
This work was prepared with the support of the "RUDN University Program 5-100".


\begin{thebibliography}{99} %


\bibitem{AU}
R. Ashurov and S. Umarov, Determination of the order of fractional 
derivative for subdiffusion equation, preprint, 
arXiv:2005.13468v1

\bibitem{CNYY}
J. Cheng, J. Nakagawa, M. Yamamoto and T. Yamazaki, Uniqueness
in an inverse problem for a one dimensional fractional diffusion equation,
Inverse Problems {\bf 25} (2009) 115002.

\bibitem{GLY}
R. Gorenflo, Y. Luchko and M. Yamamoto, 
Time-fractional diffusion equation in the fractional 
Sobolev spaces,  Fract. Calc. Appl. Anal. {\bf 18} (2015) 799-820.

\bibitem{HNWY} 
Y. Hatano, J. Nakagawa, S. Wang and M. Yamamoto, Determination of order in 
fractional diffusion equation, J. Math-for-Ind. {\bf 5A} (2013) 51-57.

\bibitem{J}
J. Janno, Determination of the order of fractional derivative and a kernel in 
an inverse problem for a generalized time fractional diffusion equation, 
Electron. J. Differential Equations {\b 28} (2016) 1-28.

\bibitem{JK} J. Janno and N. Kinash, 
Reconstruction of an order of derivative and a source term in a fractional
diffusion equation from final measurements, Inverse Problems {\bf 34(2)} 
(2018), 025007.

\bibitem{KPSV}
M. Krasnoschok, S. Pereverzyev, S.V. Siryk and N. Vasylyeva, Regularized 
reconstruction of the order in semilinear subdiffusion with memory, 
RICAM-Report 2019-05, www.ricam.oeaw.ac.at/files/reports/19/rep19-05.pdf

\bibitem{KRY}
A. Kubica, K. Ryszewska and M. Yamamoto, 
{\it A Theory of Time-fractional Partial Differential Equations},
Springer, Tokyo, 2020 (to appear).

\bibitem{KY}
A. Kubica and M. Yamamoto, 
Initial-boundary value problems for fractional diffusion
equations with time-dependent coefficients, Fract. Calc. Appl. Anal. 
{\bf 21} (2018)  276-311.

\bibitem{LZJY}
G. Li, D. Zhang, X. Jia and M. Yamamoto, Simultaneous inversion for the 
space-dependent diffusion coefficient and the fractional order in the 
time-fractional diffusion equation, Inverse Problems {\bf 29} (2013) 065014.

\bibitem{LiLiuY}
Z. Li, Y. Liu, and M. Yamamoto, Inverse problems of determining parameters
of the fractional partial differential equations, 431-442,
in Handbook of Fractional Calculus with Applications. Vol. 2
(eds. by A. N. Kochubei, Y. Luchko and J.A. Tenreiro Machado),
De Gruyter, Berlin, 2019.

\bibitem{LLY}
Z. Li, Y. Luchko and M. Yamamoto, Analyticity of solutions to a distributed
order time-fractional diffusion equation and its application to an inverse
problem, Computers and Mathematics with Applications {\bf 73} (2017) 
1041-1052.

\bibitem{LY}
Z. Li and M. Yamamoto, Uniqueness for inverse problems of determining
orders of multi-term time-fractional derivatives of diffusion equation, 
Applicable Analysis {\bf 94} (2015) 570-579.

\bibitem{LiY}
Z. Li and M. Yamamoto, Inverse problems of determining coefficients of
the fractional partial differential equations, pp.443-464,
in Handbook of Fractional Calculus with Applications. Vol. 2
(eds. by A. N. Kochubei, Y. Luchko and J.A. Tenreiro Machado),
De Gruyter, Berlin, 2019.


\bibitem{LiuLiY}
Y. Liu, Z. Li and M. Yamamoto, Inverse problems of determining 
sources of the fractional partial differential equations, 411-429,
in Handbook of Fractional Calculus with Applications. Vol. 2
(eds. by A. N. Kochubei, Y. Luchko and J.A. Tenreiro Machado),
De Gruyter, Berlin, 2019.

\bibitem{MK}
R. Metzler and J. Klafter, The random walk's guide to anomalous diffusion: 
a fractional dynamics approach. Phys. Rep. {\bf 339} No 1 (2000) 1-77.

\bibitem{Pa} A. Pazy, 
{\it Semigroups of Linear Operators and Applications to Partial 
Differential Equations}, Springer-Verlag, Berlin, 1983.

\bibitem{Po}
I. Podlubny, {\it Fractional Differential Equations}, Academic Press, 
San Diego, 1999.

\bibitem{RR}
M. Renardy and R.C. Rogers,
{\it An Introduction to Partial Differential Equations},
Springer, 1993, New York.

\bibitem{RZW}
Z. Ruan, W. Zhang and Z. Wang, Simultaneous inversion of the fractional order 
and the space-dependent source term for the time-fractional diffusion 
equation, Appl. Math. Comput. {\bf 328} (2018) 365-379.

\bibitem{RZ}
W. Rundell and Z. Zhang, Recovering an unknown source in a fractional 
diffusion problem, Journal of Computational Physics {\bf 368} (2018)
299-314.

\bibitem{SY} 
K. Sakamoto and M. Yamamoto, 
Initial value/boundary value problems for fractional diffusion-wave
equations and applications to some inverse problems, J. Math. Anal. Appl. 
{\bf 382} (2011) 426-447.

\bibitem{SLJ}
C. Sun, G. Li and X. Jia, Simultaneous inversion for the diffusion 
and source coefficients
in the multi-term TFDE, Inverse Probl. Sci. Eng. {\bf 25} (2017) 1618-1638.

\bibitem{Ta}
H. Tanabe, {\it Equations of Evolution}, Pitman, London, 1979.

\bibitem{TU}
S. Tatar and S. Ulusoy, A uniqueness result for an inverse problem in a
space-time fractional diffusion equation, Electronic Journal of Differential
Equations {\bf 258} (2013), 1-9.

\bibitem{YJQ}
B. Yu, X. Jiang and H. Qi, An inverse problem to estimate an unknown
order of a Riemann-Liouville fractional derivative for a fractional 
Stokes's first problem for a heated generalized second grade fluid, 
Acta Mech. Sin. {\bf 31} (2015) 153-161.

\bibitem{Za}
R. Zacher, Weak solutions of abstract evolutionary integro-differential
equations in Hilbert spaces, Funkcialaj Ekvacioj {\bf 52} (2009) 1-18.

\bibitem{ZCW}
X. Zheng, J. Cheng and H. Wang, Uniqueness of determining the variable 
fractional order in variable-order time-fractional diffusion equations,
Inverse Problems {\bf 35} (2019) 1-11.
\end{thebibliography}
\end{document}